\documentclass[11pt]{amsart}
\usepackage{color}
\usepackage{amssymb,verbatim}
\usepackage{hyperref}
\usepackage{amsthm,array,amssymb,amscd,amsfonts,amssymb,latexsym, url}
\usepackage{amsmath}
\usepackage[all]{xy}

\usepackage{cyrillic}

\newcommand{\cyrit}{\fontencoding{OT2}\selectfont\textcyrit}

\newcommand{\cyrbf}{\fontencoding{OT2}\selectfont\textcyrbf}

\renewcommand{\lim}{\varprojlim}

\numberwithin{equation}{section}

\newfont{\gothic}{eufb10}

\renewcommand{\qed}{{\hfill$\square$}}

\theoremstyle{plain}
\newtheorem{thm}{Theorem}[section]
\renewcommand{\proofname}{Proof}
\newtheorem{prop}[thm]{Proposition}
\newtheorem{coro}[thm]{Corollary}

\theoremstyle{definition}
\newtheorem{defi}[thm]{Definition}
\newtheorem{rem}[thm]{Remark}
\newtheorem{thmd}[thm]{Theorem}
\newtheorem{propd}[thm]{Proposition}

\newtheorem{exmp}[thm]{Example}

\newcommand{\Hn}{\ensuremath{H^2_{nr}}}
\newcommand{\nH}{\ensuremath{H_{nr}}}

\begin{document}

  \title[Stable rationality]
 {Varieties that are not stably rational, zero-cycles and unramified cohomology}

\author{Alena Pirutka}
\address{Courant Institute of Mathematical Sciences,
New York University,
251 Mercer str., New York, NY 10012}
\email{pirutka@cims.nyu.edu}

$\,$
\vspace{0.5cm}
\maketitle

\vspace{-1cm}
 \begin{abstract}
This is a survey of recent examples of varieties that are not stably rational. We review the 
specialization method based on properties of the Chow group of zero-cycles used in these examples and 
explain the point of view of unramified cohomology for the construction of nontrivial stable invariants of the special fiber. In particular, we find an explicit formula for the Brauer group of 
fourfolds fibered in quadrics of dimension $2$ over a rational surface.
      \end{abstract}

\section{Rational, unirational and stably rational varieties}

Let $X$ be a projective variety defined over a field $k$.  Several notions, describing how close $X$ 
is to projective space, have been developed:
\begin{itemize}\item $X$ is {\bf rational} if $X$ is birational to a projective space $\mathbb P^n_k$;
\item $X$ is {\bf stably rational} if $X\times \mathbb P^m_k$ is rational, for some integer $m$;
\item $X$ is {\bf retract rational} if the identity map on $X$ factorizes rationally through a projective space: 
there exist Zariski open subsets $U\subset X$ and $V\subset \mathbb P^n_k$, for some $n$, and two maps $f:U\to V$ and $g:V\to U$ such that $g\circ f=id_U$;
\item $X$ is {\bf unirational} if there is a dominant rational  map $\mathbb P^n_k\dashrightarrow X$ (one may assume $n=\mathrm{dim}\,X$, see \cite[Lemma 11]{K00});
\item $X$ is {\bf rationally connected} if for every algebraically closed field  $\Omega$ containing $k$,
any $\Omega$-points $x_1, x_2\in X(\Omega)$ can be joined by a rational curve: there exists 
a morphism $\mathbb P^1_{\Omega}\to X_{\Omega}$, with $0\mapsto x_1, \infty\mapsto x_2$. 
\end{itemize}

When $k$ is not algebraically closed, sometimes the term {\it $k$-rational}  is used, 
in order to emphasize that the birational map $X\dashrightarrow \mathbb P^n_k$ is defined over $k$ 
(similarly for the other notions but the last one).  
From the algebraic point of view, $X$ is rational if and only if the function field $k(X)$ of $X$ is a purely transcendental extension of $k$, $X$ is stably rational if and only if after adjoining some algebraically independent variables $y_1,\ldots, y_m$ the field $k(X)(y_1,\ldots, y_m)$ becomes purely transcendental over $k$ and $X$ is unirational if and only if the field $k(X)$ is a subfield of a purely transcendental extension of $k$.  Also for the first four notions one does not need to assume that $X$ is proper.
 From the definitions we deduce the following implications:
 \begin{center}
 $X$ is rational\\
 $\Downarrow$ \textcircled{\footnotesize{1}} \\
 $X$ is stably rational\\
 $\Downarrow$  \textcircled{\footnotesize{2}} \\
 $X$ is retract rational\\
 $\Downarrow$  \textcircled{\footnotesize{3}} \\
 $X$ is unirational\\
 $\Downarrow$  \textcircled{\footnotesize{4}} \\
 $X$ is rationally connected.
 \end{center}
Currently, it is still unknown if the implications   
\textcircled{\footnotesize{2}}  and  \textcircled{\footnotesize{4}} are equivalences or not, even for  $k=\mathbb C$.
It is commonly believed that they are not. Indeed, a smooth Fano variety over an algebraically closed field $k$ 
of characteristic zero is rationally connected: in particular, a hypersurface of degree $d$ in $\mathbb P^n_k$, with $d\leq n$; unirationality of all such varieties
seems unlikely. The case when $d=n$ is probably the most plausible to investigate as a 
possible counterexample to the converse of the implication \textcircled{\footnotesize{4}}.

In general, it is still unknown if a smooth hypersurface of degree $d\geq 4$ in $\mathbb P^n_{\mathbb C}$ could be  rational (resp. stably rational, resp. retract rational) \cite[p.282]{kollar}.\\

\paragraph{\bf Unirational nonrational varieties.}
Let $X$ be a smooth projective variety with $X(k)\neq\emptyset$. If $\mathrm{dim}\, X=1$, the classical L\"uroth theorem says that all the rationality notions above are equivalent to $X\simeq \mathbb P^1_k$. If $X$ is a smooth projective complex surface, it follows from the Castelnuovo-Enriques classification that all the notions above are equivalent as well. However, over nonclosed fields this is no longer true. For example, there exist nonrational smooth cubic surfaces (with a point), see \cite[V.2.7]{Manin}, and a smooth cubic hypersurface is always unirational \cite{K00}.  
It was a classical
open question to decide whether or not there exists a unirational nonrational smooth projective complex
variety of dimension at least $3$. This question was resolved in the 
$\mathrm{1970}s$, when three types of examples appeared.  
Surveys describing these examples are available (see, e.g., \cite{survolurnr}), here we  provide only a brief summary and mention what is currently known for the stable rationality of these examples:
\begin{itemize}
\item Let $X\subset \mathbb P^4_{\mathbb C}$ be a smooth cubic threefold.  As mentioned above, these varieties are unirational. Clemens and Griffiths proved that they are never rational.  To detect failure of rationality they used  the intermediate jacobian: for a rational variety it is a product  of jacobians of curves, which fails in the case of $X$. This invariant does not obstruct
stable rationality: in fact, it is still an open question whether or not a cubic threefold is stably rational. 
By recent work of B. Hassett and Y. Tschinkel \cite{HT16}, it is essentially 
the only remaining family of Fano varieties of dimension $3$ for which stable rationality is not settled (see below).
\item Let $X\subset \mathbb P^4_{\mathbb C}$ be a smooth quartic threefold. Some of these varieties are unirational, for example the quartic  $$x_0^4+x_1^4+x_2^4+x_3^4+x_0x_4^3+x_3^3x_4-6x_1^2x_2^2=0$$
(see \cite{IM} for the geometric construction going back to Segre).  
It is still unknown whether or not they are all unirational. 
As mentioned above, they are all rationally connected.
Iskovskikh and Manin \cite{IM} established that no such quartic is rational 
using the method of birational rigidity: they proved that every birational automorphism of $X$ extends to an 
automorphism of $X$. However, the group of birational automorphisms of rational varieties is huge. The birational rigidity method has been extensively studied and applied, see in particular the works of Pukhlikov (and his survey 
 \cite{Pukhlikov}) and Cheltsov (and the survey \cite{Cheltsov}). 
Recently, the method has been further developed by T. de Fernex, 
who proved that no smooth hypersurface of degree $n\ge 4$ in $\mathbb P^n_{\mathbb C}$ is rational \cite{dF}.

Using the specialization method we discuss in this survey, in joint work with 
J.-L. Colliot-Th\'el\`ene we showed that a very general quartic threefold is not stably rational \cite{CTP15}.

\item Artin and Mumford \cite{artinmumford} constructed a unirational variety  that is not even stably rational: one could view their examples either as a double cover $$X: z_4^2-f(z_0, z_1, z_2, z_3)=0$$ of $\mathbb P^3_{\mathbb C}$ ramified along a (particular) quartic $f(z_0, z_1, z_2, z_3)=0$
 or as a conic bundle over $\mathbb P^2_{\mathbb C}$. The nontrivial invariant they use is the $2$-torsion in the Brauer group 
$Br(X)[2]$ (which in this case is isomorphic to $H^3(X,\mathbb Z)[2]$ as well). For stably rational varieties this invariant vanishes. This example is one of the key ingredients of proofs of failure of stable rationality in the classes
of varieties we discuss  in this paper.\\
\end{itemize}

\paragraph{\bf Stably rational nonrational varieties.} Examples of stably rational nonrational complex varieties $X$  were constructed in joint work of Beauville,  Colliot-Th\'el\`ene,  Sansuc and  Swinnerton-Dyer  \cite{BCTSSD}. 
These examples are also threefolds admitting a conic bundle structure over a rational surface. 
The method to establish that $X$ is nonrational relies on the intermediate jacobian criterion. 
Showing that $X$  is stably rational uses the universal torsors techniques \cite{CTSan79}.   
If $k$ is not algebraically closed, there are examples of stably rational nonrational  surfaces
as soon as $k$ admits a Galois extension $K$ with $Gal(K/k)\simeq \mathfrak{S}_3$  \cite{BCTSSD}. 
However, the following question of B. Hassett remains open:  

\vspace{0.1cm}
{ \it is there a stably rational nonrational surface $S$ defined over a finite field $k$}?\\

\paragraph{\bf Algebraic groups and classifying spaces.}
Let $G$ be a linear algebraic group over $k$. There is a generically free representation $V$ of $G$ and an open $G$-equivariant subset $U\subset V$ such that there is a $G$-torsor $U\to U/G$; in addition, if $U'\subset V'$ is another such open in a generically free representation of $G$, then $U/G$ and $U'/G'$ are stably birational \cite{CTSa07}.  One then views $U/G$ as an algebraic approximation of the classifying space $BG$ of $G$ and one says that $BG$ is stably (resp. retract) rational if $U/G$ is \cite{Merkurjev}. 
Until now, no example of non retract rational $BG$, for $G$ a connected linear algebraic group over an algebraically closed field, 
has been found. A plausible candidate could be $G=\mathrm{Spin}_{16}$ (see  \cite{Merkurjev}). However, such examples exist  if $G$ is finite (Saltman \cite{Saltman1}, Bogomolov \cite{Bogomolov}) or if $k$ is not algebraically closed (Merkurjev \cite{Merkurjev1}).  For $G=\mathrm{PGL}_p$ with $p$ a prime,  the classifying space $BG$ is retract rational \cite{Saltman}. It is an open question whether or not 
it is stably rational. 

\vspace{0.7 cm}

This survey is devoted to recent constructions of examples of varieties that are not stably rational, complementing      
surveys of Arnaud Beauville \cite{BeauvilleRapport} and Claire Voisin \cite{VoisinRapport} on the same topic. 
In section \ref{nonstr} we provide a list of new examples and 
explain the specialization method that allowed to produce them.  
The example of Artin and Mumford and their computations of the Brauer group of fibrations in conics  
are crucial steps in the construction of many cases. 
We review a method to compute the Brauer group as an unramified cohomology group 
and discuss some applications to fibrations in quadrics. 
In particular, we give an explicit formula for the Brauer group of 
fourfolds fibered in quadrics of dimension $2$ over a rational surface.\\
 
 \paragraph{\bf Acknowledgements.}  I benefitted very much from motivating discussions with 
 Jean-Louis Colliot-Th\'el\`ene, Brendan Hassett, Yuri Tschinkel  and Asher Auel. I am very grateful to
 Alexander Merkurjev for carefully reading a preliminary version of this paper and for correcting several inaccuracies.

 \vspace{1cm}

\section{Specialization method and applications}\label{nonstr}

\subsection{Examples}\label{sectionEx}

We now present recent examples of varieties failing stable rationality, in chronological order of their discovery. 
These varieties are not even retract rational, but the statements are often in terms of the more intuitive property of stable (non) rationality.   
{\it Very general} means that we remove a countable set of proper 
closed conditions on coefficients of the corresponding varieties\footnote{More precisely, if  $B$ is an integral  scheme over an algebraically closed field $k$ and $\mathcal X\to B$  is a family of algebraic varieties over $B$, we say that a very general member of the family is not stably rational if there exists a countable set  of proper closed subvarieties $\{B_i\subset B\}_{i\in I}$ such that for any $b\in B(k)\setminus (\bigcup_{i\in I}B_i(k))$, the variety $\mathcal X_{b}$ is not stably rational.}. 
\begin{enumerate}
\item  (C. Voisin \cite{V13}) A very general  double cover $X\to \mathbb P^3_{\mathbb C}$, ramified along a quartic with at most $7$ nodes.
\item (Colliot-Th\'el\`ene -- Pirutka \cite{CTP15}) A very general quartic threefold $X_4\subset \mathbb P^4_{\mathbb C}$.
\item (Beauville  \cite{B6})  A very general  double cover $X$  of $\mathbb P^3_{\mathbb C}$ branched along a  sextic.  
\item  (Totaro \cite{To15}) A very general hypersurface $X_d\subset \mathbb P^{n+1}_{\mathbb C}$, $n\geq 3$, of degree  $d\geq 2\lceil (n+2)/3 \rceil$. 
\item (Hassett--Kresch--Tschinkel  \cite{HKT15})  A very general conic bundle over  $\mathbb P^2_k$ with $k=\bar k, char\,k\neq 2$, with  discriminant curve of degree $d\geq 6$.  More generally, let $S$ be a smooth projective rational surface over $k$.  Using \cite[Thm 1.1]{artinmumford} (see also section \ref{fq} below), the data of a conic bundle $X$ over $S$ with at most  nodal (reduced) discriminant curve $D\subset S$ is equivalent to the data of $D$ and a degree $2$ \'etale covering of $D$.   The space of such curves in a fixed linear system  $\mathcal L$ of effective divisors on $S$, together with  a degree $2$ \'etale covering, is not necessarily irreducible. Let $\mathcal M$ be an irreducible component. Assume that a general member of $\mathcal L$ is irreducible and that $\mathcal M$ contains a reducible curve $D$ 
with smooth irreducible components, with a cover that is nontrivial on each  component. Then the conic bundle corresponding to a very general point of $\mathcal M$ is not stably rational.
 \item   (Colliot-Th\'el\`ene -- Pirutka \cite{CTP15-2}) A very general cyclic cover $X$ of $\mathbb P^n$  of prime degree $p$ ramified along a  hypersurface of degree $mp$ with $m(p-1)<n+1\leq mp$. Here the condition $m(p-1)<n+1$ is just to insure that $X$ is rationally connected, and the condition $n+1\leq mp$  insures that $X$ is not stably rational.
\item   (Hassett--Tschinkel  \cite{HT16})  A very general nonrational smooth Fano threefold $X$, if $X$ is not birational to a cubic threefold. Here very general refers to families of such Fano threefolds with some fixed numerical invariants, such as "hypersurface in $\mathbb P(1,1,1,2,3)$ of degree $6$", "quartic in $\mathbb P^4$\," etc., see \cite[Section 2]{HT16} for more details.
\end{enumerate}

The  stable birational invariants, essentially of cohomological nature,  that have been successfully used in other contexts previously, are either trivial  in many of these examples, or are difficult to compute (see remarks \ref{trhyp} or \ref{trirr} below).  The strategy in the specialization method is as follows:
\begin{itemize}
\item[step $1$:] use the invariants that behave well under specialization,  even to mildly singular varieties;
\item [step $2$:] produce some varieties (special fibers) with nontrivial invariants;
\item [step $3$:] show that the varieties of interest specialize to the examples produced in the previous step.\\
\end{itemize}

For example, the invariant $H^3(X, \mathbb Z)[2]$ of Artin and Mumford  does not behave well under specialization:  it is trivial for a general quartic double solid. In the next section, we introduce invariants based on the properties of zero-cycles,  satisfying the requirements of the step $1$ above.

\vspace{1cm}

\subsection{Zero-cycles and universal triviality of $CH_0$}$\,$\\

Let $X$ be a projective variety of dimension $n$ defined over a field $k$.   The decomposition properties of the class $N[\Delta_X]$,  $N>0$,    of the diagonal $\Delta_X=\{x,x\}\subset X\times X$, in the Chow group $CH_n(X\times X)$, have been previously used for various applications (see for example  the work of Bloch and  Srinivas \cite{BlochSrinivas}).    Here we consider   the following integral version:

\defi{  A projective variety $X$ of dimension $n$ defined over a field $k$, with $X(k)\neq \emptyset$, has {\bf Chow decomposition of the diagonal}, if 
\begin{equation}\label{diagdec}
[\Delta_X]=[X\times x]+[Z]\in CH_n(X\times X),
\end{equation} where $x\in X(k)$ and $Z$ is an $n$-cycle on $X\times X$, supported on $D\times X$, where $D\subset X$  is a closed subvariety of codimension  at least $1$. }

\rem{One could replace $x\in X(k)$ by a zero-cycle of degree $1$. If $X$ is smooth, the existence of the decomposition above does not depend on the choice of $x$ (see \cite[Lemma 1.3]{CTAP}).  Also  assuming $X$ proper is enough. \\}

The following notion is more flexible for the properties of specialization and for varieties with mild singularities (see also \cite{Merkurjev2}).

\defi{ A projective variety $X$ is {\bf universally $CH_0$-trivial} if for any field  $F$ containing $k$ the degree map $CH_0(X_F)\to \mathbb Z$ is an isomorphism.\\}

We also have a relative notion:

\defi{A projective morphism $f:X\to Y$ of varieties over $k$ is {\bf universally $CH_0$-trivial} if for any field $F$ containing $k$ the push-forward map $f_*: CH_0(X_F)\to CH_0(Y_F)$ is an isomorphism. \\}

\exmp{\begin{enumerate}
\item If $X$ is a smooth, projective, retract rational variety, then $X$ is universally $CH_0$-trivial \cite[Lemma 1.5]{CTP15}.
\item  Some singular or reducible examples of  universally $CH_0$-trivial varieties:
\begin{enumerate}
\item  a connected variety $X$ over an algebraically closed field, such that each component of  the reduced scheme $X^{red}$ is a rational variety with isolated singular points \cite[Lemma 2.3]{CTP15-2};
\item a projective, reduced, geometrically connected variety $X=\cup_{i=1}^n X_i$, such that each $X_i$ is   universally $CH_0$-trivial and geometrically irreducible, and each intersection $X_i\cap X_j$ is either  empty or has a $0$-cycle of degree $1$ \cite[Lemma 2.4]{CTP15-2}.
\end{enumerate}
\item If $f:Z\to Y$ is a projective morphism, such that for any (scheme) point $M$ of $Y$, with residue field $\kappa(M)$, the fiber $Z_{\kappa(M)}$ is  universally $CH_0$-trivial, then $f$ is universally $CH_0$-trivial  \cite[Prop. 1.7]{CTP15} (note that the converse  property does not hold). This criterion is useful in particular when $f$ is the resolution of singularities  of $Y$, and the exceptional divisors are of type described in examples above. In particular, this applies if $Y$ has just ordinary double point singularities.
\item Let $X$ be a projective integral variety over a field $k$. If $X$ is universally $CH_0$-trivial, then $X$ admits  Chow decomposition of the diagonal: for $F=k(X)$, $\eta\in X(k(X))$  the generic point and  $x\in CH_0(X)$  a zero-cycle of degree $1$, we have that the zero-cycle $\eta-x\in CH_0(X_F)$   is of degree $0$, hence $\eta-x=0$. This means that the restriction of the cycle $[\Delta_X]-[X\times x]$ to $CH_0(X_{k(X)})$ (with respect to the first projection) is zero, hence this cycle is supported on $D\times X$ as in the definition. If in addition $X$ is smooth, then these notions are equivalent (see below).  
\item  If $X$ is a rationally connected variety over a field $k$, then  $X$ is not necessarily universally $CH_0$-trivial. Indeed, we only have that for any algebraically closed field $F$ containing $k$ the degree map $CH_0(X_F)\to \mathbb Z$ is an isomorphism.  For  rigidity properties of the Chow groups for varieties over algebraically closed fields see \cite{Le}. \\
\end{enumerate}
}

The point of view of decomposition of the diagonal is very useful to establish the triviality of various invariants: if $X$ is smooth and projective,  we have that the group $CH_n(X\times X)$ acts on many invariants of $X$ (such as $H^3(X, \mathbb Z)_{tors}$, $Br(X)$ etc.), the diagonal $[\Delta_X]$ acts as identity.  The decomposition (\ref{diagdec}) often allows to conclude that the action $[\Delta_X]_*$ is a zero map, so that the invariants above are trivial:

\thm\label{intdectrinv}{Let $X$ be a smooth projective variety, defined over a field $k$. Assume that $X$ has Chow decomposition of the diagonal (\ref{diagdec}). Then \begin{itemize}
\item[(i)]  $X$ is universally $CH_0$-trivial;
\item[(ii)] for any field $F$ containing $k$, the natural map $Br(F)\to Br(X_F)$ is an isomorphism;
\item[(iii)]  more generally, for any Rost cycle module $M^i$ over $k$ and for any field $F$ containing $k$, the natural map $M^i(F)\to M^i_{nr}(F(X)/F)$ is an isomorphism.
\item[(iv)]  $H^0(X, \Omega_X^i)=0$ for any $i>0$.
\end{itemize} }
\proof{For $(i)-(iii)$ see \cite[Section 1]{CTP15} (and references there). For $(iv)$, see \cite[Lemma 2.2]{To15}, and the remark below.\qed\\} 

\rem{\begin{enumerate}
\item For the definition and the properties of the Rost cycle modules, see \cite{Rost}.  In particular, the unramified cohomology groups, defined in \ref{defnr} below, are the groups $M^i_{nr}$ for the cycle module of the Galois cohomology.
\item  If $k=\mathbb C$, the property $(ii)$ also implies that $H^3(X, \mathbb Z)_{tors}=0$.
\item In order to construct the action  of correspondences on  $H^j(X, \Omega_X^i)$ one may use cycle class maps $CH^i(X)\to H^i(X, \Omega_X^i)$  as in \cite{EZ}, \cite[Section 2]{Srinivas}.
\item More precisely, B. Totaro proves that if the group  $H^0(X, \Omega^i_X)$ is nontrivial for some $i>0$, then the group $CH_0(X)\otimes \mathbb Q$ is not universally trivial if $k$ is of characteristic zero, and the group  $CH_0(X)/p$ is not universally trivial if $k$ is of positive characteristic $p$.\\
\end{enumerate}}

\subsection{Specialization properties}
The point of view of zero-cycles allows to establish the following local specialization property (see \cite[Thm. 1.11, Thm. 1.12]{CTP15}):

\thmd\label{localsp}{\it Let $A$ be a discrete valuation ring, $K$ the field of fractions of $A$ and $k$ the residue field. Let $\mathfrak X$ be a faithfully flat, proper $A$-scheme, with integral geometric fibers. Let $X/K$ be the generic fiber of $\mathfrak X$ and $Y/k$ the special fiber. Assume that there is a resolution of singularities $f:Z\to Y$ of $Y$ (i.e. $f$ is a proper birational map and $Z$ is smooth), such that $f$ is universally $CH_0$-trivial. Assume one of the following:
\begin{itemize}
\item[(i)] $X$ is smooth, universally $CH_0$-trivial and $Z$ has a zero-cycle of degree $1$;
\item[(ii)] $k$ is algebraically closed and for $\bar K$ an algebraic closure of $K$, the variety $X_{\bar K}$ is universally $CH_0$-trivial.
\end{itemize}
Then $Z$ is universally $CH_0$-trivial.\\
}

\rem\label{remappli}{\begin{enumerate}
\item The proof of the theorem above uses the existence of the specialization map $CH_0(X)\to CH_0(Y)$ for the Chow group of zero-cycles \cite{fulton}. In general the specialization maps for Chow groups are subtle to define, but  here  $Y$ is a Cartier divisor on $\mathfrak X$ and  \cite[Prop. 2.6]{fulton} provides a simple construction.
\item  In (i) (resp. (ii))  one could replace the assumption that $X$ is smooth by ``$X$ (resp. $X_{\bar K}$) has a resolution $\tilde X\to X$ with $\tilde X$ smooth, projective,  and universally $CH_0$-trivial. "
\item Burt Totaro  \cite{To15} has generalized the theorem above to the case when $Y$ is not necessarily irreducible: under assumption (ii), he shows that for any extension $L$ of $k$, every zero-cycle of degree $0$ in the smooth locus of $Y_L$ is zero in $CH_0(Y_L)$.
\item The applications of the theorem are usually as follows: assume that $Y/k$ is a singular variety, such that $Y$ has a resolution $f:Z\to Y$ with $f$ universally $CH_0$-trivial.  Assume that we know that for any smooth model $Z$ of $Y$ some birational invariant (for example, $Br(Z)$) is nontrivial. Then every variety $X$ that specializes to $Y$ (i.e. $X$ and $Y$ are fibers of a local family $\mathfrak X$ as above), is not retract rational.\\
\end{enumerate}}

The property of the diagonal decomposition is more suitable for the following global property (see \cite[Thm.2.1]{V13}, \cite[Thm. 2.3]{CTP15}):

\thmd\label{globalsp}{\it Let $B$ be an integral scheme of finite type over an uncountable algebraically closed field $k$ of characteristic zero. Let $f:\mathfrak X\to B$ be a flat projective map. If there is a point $b_0\in B(k)$ such that $\mathfrak X_{b_0}$ has no Chow decomposition of the diagonal, then for a very general point  $b\in B(k)$ the fiber $\mathfrak X_{b}$ has no Chow decomposition of the diagonal.\\ }

\rem\label{remglob}{\begin{enumerate}
\item In particular, if there is a smooth fiber $\mathfrak X_b$ of $f$, such that $\mathfrak X_b$ is not universally $CH_0$-trivial, then a very general fiber of $f$ is not universally $CH_0$-trivial.
\item   The theorem above could also be viewed  as  a property of the geometric generic fiber of $f$: in fact, by a general statement  \cite[Lemma 2.1]{Vial},  a  very general fiber of $f$ is isomorphic to the geometric generic fiber of $f$, as an abstract scheme.
\item By \cite[Thm. 9]{HKT15}, the theorem also holds in characteristic $p$ if one considers the decomposition of the diagonal with $\mathbb Z[\frac{1}{p}]$-coefficients.
\item It should be possible to replace the condition that $f$ is flat by a weaker condition, related to the existence of a relative cycle  (in the sense of Koll\'ar \cite{kollar}) of the diagonal $\Delta_{\mathfrak {X}}$ over $B$.\\
\end{enumerate}
} 
 
 \subsection{Point of view of $R$-equivalence.}
 
Let $X$ be a projective variety over a field $k$. Recall that two points $x,y\in X(k)$ are {\bf directly $R$-equivalent} if there is a map $f:\mathbb P^1_k\to X$, $0\mapsto x, \infty\mapsto y$. This generates an equivalence relation called {\bf $R$-equivalence}.  The set of $R$-equivalence classes is denoted $X(k)/R$. 
 
 An analogue of the local specialization property \ref{localsp} could also be obtained using the (a priori, more ``elementary" to define) set $X(k)/R$ instead of the Chow group of zero-cycles. Theorem \ref{principalR} below suffices for some applications.

\begin{prop}\label{specialisationR}

Let $A$ be a discrete valuation ring, $K$ the field of fractions of $A$ and $k$ the residue field. Let $\mathfrak X$ be a proper $A$-scheme. Let $X/K$ be the generic fiber of $\mathfrak X$ and $Y/k$ the special fiber. \begin{itemize}
\item [(i)] The natural map
$X(K)={\mathfrak X}(A) \to Y(k)$ induces the map
$$
X(K)/R \to Y(k)/R.
$$
\item [(ii)] If ${\mathfrak X}/A$ is proper and flat, with  $A$ henselian,
and $X(K)/R$ has at most one element
then the image of  $Y^{sm}(k)$ in $Y(k)/R$ has at most one element.
Here $Y^{sm}\subset Y$ denotes the smooth locus of $Y$. 
\end{itemize}
\end{prop}
\proof{
For (i) see \cite{madore}. The statement (ii) follows from (i). \qed \\
}

The following condition replaces the existence of a $CH_0$-trivial resolution (see \cite[Rem. 1.19]{CTP15}):
\begin{prop}\label{bonneresolution}
 Let $f : Z \to Y$ be a  birational morphism of  integral proper  $k$-schemes. Assume that there exists an open
 $U \subset Y$ such that 
$V= f^{-1}(U)$ satisfies the following properties:
\begin{itemize}
\item[(i)] The induced map
$f : V \to U$ is an isomorphism.

\item[(ii)]  The complement of $U$  in $Y$ is a finite union of geometrically integral $k$-schemes $F_{i}$,
such that for each $S_{i}=f^{-1}(F_{i}) \subset Z$, the set $S_{i}(k)/R$ has at most one element.

\item[(iii)]  The image of  $U(k)$ in $Y(k)/R$ has at most one element.
\end{itemize}
Then the image of $V(k)$ in $Z(k)/R$  has at most one element. \qed \\
\end{prop}

\thm{\label{principalR}
Let $A$ be a discrete valuation ring, $K$ the field of fractions of $A$ and $k$ the residue field. Assume $k$ is algebraically closed of characteristic zero.
Let $\mathfrak X$ be a faithfully flat proper $A$-scheme.
Assume that
\begin{itemize}\item[(i)] the generic fiber $X=\mathfrak X\times_{A}K$  is  geometrically integral, smooth and  $X_{\bar K}$ is stably rational;
\item[(ii)] the special fiber $Y={\mathfrak X}\times_{A}k$ is integral, and there exist a nonempty smooth open  $U \subset Y$ and a resolution of singularities $f : Z \to Y$ such that for the induced map $f : V:=f^{-1}(U) \to U$
the properties  (i) and (ii) of Proposition \ref{bonneresolution} are satisfied over any field 
  $L$ containing  $k$.
\end{itemize}

Then $Br(Z)=0$. 
More generally, for any Rost cycle module $M^i$ over $k$ and for any field $F$ containing $k$ and for $i \geq 1$, we have $M^{i}_{nr}(F(Z_F))=0$.
}

\proof{
We may assume that $A=k[[t]]$
and that there exists a finite extension $K'/K$ such that $X_{K'}$ is stably rational (over  $K'$).
The integral closure $A'$ of $A$ in $K'$ is a discrete valuation ring
with residue field $k$. 
Let $L$  be the field of functions of $Z$. The generic point of $Y$ defines an $L$-point $\eta$ of
 $Z_{L}$, this point is in fact in
 $V_{L}(L)$. Pick a $k$-point of $U$ and let $P$  be the corresponding $L$-point of $Z_{L}$, this point is again in $V_{L}(L)$.
Let $\alpha \in M^{i}_{nr}(k(Z))$.
By base change from $k[[t]]$ to $L[[t]]$, the assumptions of Proposition \ref{specialisationR} are satisfied. We deduce that the image of $V(L)$ in $Y(L)/R$ has at most one point.
From (ii), the $L$-points $\eta$ and $P$ of $V_{L}(L)$ are $R$-equivalent over
$Z_{L}$.  Restricting to a rational curve passing through $\eta$ and $P$ and using the triviality of $M^{i}_{nr} $ of $\mathbb P^1_k$,  we  deduce that
 $\alpha(\eta) = \alpha(P) \in M^{i}(L)$.
 But $\alpha(P)$ comes from $M^{i}(k)=0$, hence 
 $\alpha \in M^{i}(k(Z))$ is zero as well.
 \qed\\
}

\subsection{Varieties with nontrivial invariants}$\,$\\\label{badvar}

\noindent {\it Example of Artin and Mumford.\\}
Let $\mathcal D\subset \mathbb P^2_{\mathbb C}$ be a smooth conic defined by a homogeneous equation $\delta(z_0,z_1,z_2)=0$.
Let $E_1, E_2\subset  \mathbb P^2_{\mathbb C}$ be smooth elliptic curves defined by $\epsilon_1(z_0,z_1, z_2)=0$ and $\epsilon_2(z_0, z_1, z_2)=0$, respectively, each tangent to $\mathcal D$ at three points, all distinct. Assume that the curves 
$E_{1}$ and  $E_{2}$ have  nine pairwise distinct intersection points,  and distinct from the previous ones.  Then there are homogeneous forms $\beta(z_0, z_1, z_2)$ and $\gamma(z_0, z_1, z_2)$ of degrees respectively $3$ and $4$, such that $$\beta^2-4\delta\gamma=\epsilon_1\epsilon_2.$$Let $V$ be a double cover of $\mathbb P^3_{\mathbb C}$, defined by
\begin{equation}
z_4^2-\delta(z_0,z_1,z_2)z_3^2-\beta(z_0,z_1,z_2)z_3-\gamma(z_0, z_1, z_2)=0.
\end{equation}
Note that with this construction we could have $V$ defined over $\bar {\mathbb Q}$.
This double cover $V$ has $10$ ordinary double point singularities. By projection to the plane with coordinates $[z_0:z_1:z_2],$  the variety $V$ is birational to a conic bundle over $\mathbb P^2_{\mathbb C}$. 
  Artin and Mumford \cite{artinmumford} proved that for any smooth and projective variety $Z$, birational to $V$, one has $Br(Z)[2]\neq 0$. See also \cite{CTOj} for a way to prove  it using the techniques of Section \ref{fq}, as well as other examples.\\

{\it Cyclic covers.}
Let $p$ be a prime and $f(x_0, \ldots, x_n)$  a homogeneous polynomial of degree $mp$ with coefficients in an algebraically closed  field $k$.
{\bf A cyclic cover of } $\mathbb P^n_k$, branched along  $f(x_0, \ldots, x_n)=0$, is a subvariety of $\mathbb P(m,1,1,\ldots ,1)$ given by $$y^p-f(x_0, \ldots , x_n)=0.$$

If in addition the characteristic of $k$ is $p$, these covers have the following nontrivial invariants: by
  \cite[V.5.7, V.5.11]{kollar}, if
 $q: Y\to \mathbb P^n_k$ is a cyclic, degree $p$ cover of $\mathbb P^n_k$, $n\geq 3$,  branched along a hypersurface $f=0$, then for a general choice of coefficients of $f$, there  is a resolution of singularities $\pi: Z\to Y$ of $Y$ obtained by successive blow-ups of singular points,    with  $\pi^*q^*\mathcal O_{\mathbb P^n}(mp-n-1)$  a subsheaf of  $\Omega_{Z}^{n-1}$.  
In particular, 
 if   $mp-n-1\geq 0$, then
$H^0(Z, \Omega_{Z}^{n-1})\neq 0$. 
By \cite[Thm. 3.7]{CTP15-2}, one may also assume that $\pi$ is a universally $CH_0$-trivial map.
\\

More generally,  a similar analysis applies to covers of hypersurfaces (see \cite{kollar}). If, in addition, we restrict to the case $p=2$, then the singularities are just ordinary double points. So let $k$ be an algebraically closed field of characteristic $2$ and let $Y$ be a double cover of a hypersurface $\{g=0\}\subset \mathbb P^{n+1}_k$, of degree $m$, branched along $f=0$ with $deg(f)=2m$. Assume  $m\geq \lceil\frac{n+2}{3}\rceil$  and $n\geq 3$.  Then for a general choice of coefficients of $f$, there is resolution of singularities $\pi: Z\to Y$ of $Y$ obtained by successive blow-ups of (ordinary double) singular points,   such that
$H^0(Z, \Omega_{Z}^{n-1})\neq 0$. 

In Koll\'ar's original results, the assumption is slightly stronger:  $m\geq \lceil\frac{n+3}{3}\rceil$. This last condition implies that there is a big line bundle inside $\Omega_{Z}^{n-1}$, and hence $Z$ is not separably uniruled. Here we just need that $H^0(Z, \Omega_{Z}^{n-1})$ is nonzero.\\

\subsection{Back to examples}$\,$\\
In this section we give two  applications of the  technique above, in the case of hypersurfaces. \\

\noindent {\it Quartic threefolds} (see \cite{CTP15}).\\ Let $X\subset \mathbb P^4_{\mathbb C}$ be a quartic hypersurface. For a very general choice of coefficients, $X$ is not stably rational (not even retract rational):
\begin{enumerate}
\item The quartic $$Y: z_0^2z_4^2-\delta(z_0,z_1,z_2)z_3^2-\beta(z_0,z_1,z_2)z_3-\gamma(z_0, z_1, z_2)=0$$
is birational to the variety $V$ of Artin and Mumford, defined as above (by taking $z_0=1$.)  Note that $Y$ could be defined over $\bar {\mathbb Q}$. Hence, for any resolution $f:Z\to Y$ ($Z$ is smooth and projective,  $f$ is birational), one has $Br(Z)[2]\neq 0$. The singularities of $Y$ are slightly more complicated than ordinary double points, but by \cite[Appendice A]{CTP15}, there is   an explicit resolution $f$, such that $f$ is a universally $CH_0$-trivial morphism.  Moreover,   the conditions $(i)-(iii)$ of Proposition \ref{bonneresolution} are also satisfied.
\item By \ref{remappli}(3), any smooth variety that specializes to $Y$ cannot be stably rational.  For example, consider  
the projective space $\mathcal P$ corresponding to the coefficients of quartics $X\subset \mathbb P^4_{\mathbb C}$, and let $\mathcal X\to \mathcal P$ denote the universal family of quartics,  $\mathcal U\subset \mathcal P$ the open corresponding to smooth quartics and $M\in \mathcal W:=\mathcal P\setminus \mathcal U$ the $\bar{\mathbb Q}$-point corresponding to $Y$. Consider a line $L$ in $\mathcal P$ containing $M$ and not contained in $\mathcal W$.  Applying the local specialization theorem  (Theorem \ref{localsp}) to the local ring $A=\mathcal O_{L, M}$ of $L$ at $M$ we see that the smooth quartic $\mathcal X_{\bar \eta}=\mathcal X_{\overline{\mathbb C(t)}}$ corresponding to the geometric generic fiber of $L$, is not stably rational. But then for any point $P$ of $L(\mathbb C)$ that is not defined over $\bar{\mathbb Q}$, the corresponding quartic $\mathcal X_P$ is isomorphic to $\mathcal X_{\bar \eta}$, as an abstract scheme (there is a 'transcendental' parameter $P$ in the coefficients of $\mathcal X_P$, that maps to $t$ under this isomorphism). In this way we obtain 'by hand' many smooth quartics that are not stably rational. 
\item At the next step,  application of the global specialization theorem \ref{globalsp}, as explained in Remark \ref{remglob}(1) then shows that a very general quartic threefold is not stably rational. One could also show directly that $Y$ (even though $Y$ is singular) has no Chow decomposition of the diagonal (if not, $Z$ would have one with an "error" term, and it is still enough to deduce in this case that $Br(Z)=0$), see \cite[Prop. 2.4]{CTP15} and then deduce, by Theorem \ref{globalsp} that a very general quartic threefold is not stably rational.
\item Let $B\subset \mathcal U$ be the locus corresponding to quartics that are not stably rational. In  (2) we  constructed many concrete points in $B$, showing in particular that $B$ is Zariski dense in $\mathcal U$. By (3), a very general point is in $B$, but no specific point can be checked by arguments in (3). It is  a curious fact that by local and global arguments we do not obtain the same description of $B$, a priori. 
\item By a more careful specialization, we can show that $B(\bar{\mathbb Q})$ is nonempty. In fact, we have $Br(Z)=H^3_{\acute{e}t}(Z, \mathbb Z_2)[2]\neq 0$. Since $Y$ and $Z$ are defined over $\bar{\mathbb Q}$, we can extend to a family $g:\mathcal Z\to \mathcal Y$ over the ring of integers $\mathcal O_{K}$ of some number field $K$. Over an open $S\subset \mathcal O_K$, the fibers of $g$ are of the same type as the fibers of $f$, in particular, universally $CH_0$-trivial, and $H^3_{\acute{e}t}(\mathcal Z_s, \mathbb Z_2)[2]\neq 0$ for $s\in S$. Hence any smooth quartic (over $\bar{\mathbb Q }$) specializing to one of these fibers for $s\in S$ is not stably rational.
\item Choosing another specialization than to $Y$, one can also show that $B(\mathbb Q)$ is not empty.  See below.\\

\end{enumerate}

\noindent {\it Hypersurfaces in $\mathbb P^{n+1}$ of degree  $d\geq 2\lceil (n+2)/3 \rceil$ }(see \cite{To15}). 

\begin{itemize}
\item[$d=2m$:] Consider the example as in Section \ref{badvar}:
$k$ is an algebraically closed field of characteristic $2$, $Y$ is a double cover of a hypersurface $\{g=0\}\subset \mathbb P^{n+1}_k$, of degree $m$, branched along $f=0$ with $deg(f)=2m$. The condition  $m\geq \lceil\frac{n+2}{3}\rceil$  is satisfied, so that choosing general coefficients, we may assume that  there is a resolution of singularities $\pi: Z\to Y$ of $Y$ with $H^0(Z, \Omega_{Z}^{n-1})\neq 0$. Note that the map  $\pi$ is universally $CH_0$-trivial.  Let $A$ be a discrete valuation ring with  field of fractions $K$ of characteristic zero and with residue field $k$ and $t$ the uniformizing parameter of $A$. Let $F,G$ be polynomials with coefficients in $A$ that specialize to $f, g$ respectively.  Let $\mathfrak X\to Spec\,A$ be the complete intersection defined by  $y^2=F, G=ty$.  The generic fiber is a hypersurface $X$ defined by $F^2-t^2G^2=0$,   of degree $d$ and the special fiber is the double cover $Y$. By \ref{intdectrinv} and \ref{localsp}, $X_{\bar K}$ is not universally  $CH_0$-trivial. Choosing an embedding $\bar K\subset \mathbb C$,  we deduce that a very general complex hypersurface of degree $d$ is not universally  $CH_0$-trivial, hence not stably rational, by Theorem \ref{globalsp}.
\item[$d$ odd:]  Let $X$ be as in the previous case and $H\subset \mathbb P^{n+1}$ a general hyperplane.  Since one could degenerate a smooth hypersurface of degree $d$ to a union $X\cup H$, using \ref{remappli}(2) it suffices to show that $CH_0((X\cap H)_F)\to CH_0(X_{F})$ is not surjective, for some extension $F$ of the base field, for instance, $F$ the function field of $X$. Put $\bar X=X\cap H$ and let $\bar Y\subset Y$ be the special fiber of the closure of $\bar X$ in $\mathfrak X$. Choosing $H$ generally, we may assume that $\bar Y$ is  normal, and disjoint from the singular locus of $Y$. Note that $\bar Y$ is also a double over of a hypersurface. In particular, since $Y$ (resp. $\bar Y$) are normal, one could define the canonical sheaves  $K_Y$ (resp. $K_{\bar Y}$) and   $K_Y=\mathcal O(-n-2+2m)$ and  $K_{\bar Y}=\mathcal O(-n-1+2m)$, in particular $H^0(\bar Y, K_{\bar Y})=0$.  If $\bar Z$ is a resolution of singularities of $\bar Y$, that is an isomorphism over the smooth locus $\bar Y^{sm}$ of $\bar Y$, then $H^0(\bar Z, K_{\bar Z})=0$ as well: in fact, any nonzero section $\sigma\in H^0(\bar Z, K_{\bar Z})$  restricts to a nonzero section on $\bar Y^{sm}$ and any such section extends to $K_{\bar Y}$ since $\bar Y$ is normal. Now, if $CH_0(\bar X_F)\to CH_0(X_{F})$ is surjective, then one deduces that the diagonal of $X$  (resp. $Z$) has a decomposition, which is weaker than the integral diagonal decomposition we consider, but still enough to use the action of correspondences to show that the identity action on the nonzero group $H^0( Z, \Omega_{Z}^{n-1})$ factorizes through $H^0( \bar Z, \Omega_{\bar Z}^{n-1})=H^0(\bar Z, K_{\bar Z})=0$ (see \cite{To15}) and we obtain a contradiction: by \ref{remappli}(2), any smooth hypersurface degenerating to $X\cup H$ is not universally $CH_0$-trivial.
\end{itemize}

In particular, as explained in the work of B.Totaro \cite{To15}, the following quartic threefold over $\mathbb Q$:
$$(x_0x_1+x_2x_3+x_4^2)^2-4(x_0^3x_2+x_1^3x_2+x_0x_1^2x_4+x_0x_2^2x_4+x_3^3x_4)-8x_2^4=0$$
is smooth and is not stably rational, since it degenerates to the double over $f=0$ of the hypersurface $g=0$ over $\mathbb F_2$ with
$$f=x_0^3x_2+x_1^3x_2+x_0x_1^2x_4+x_0x_2^2x_4+x_3^3x_4, g=x_0x_1+x_2x_3+x_4^2.$$ This double cover has a universally $CH_0$-trivial resolution $Z$ with $H^0(Z, \Omega_{Z}^{n-1})\neq 0$. 

\rem{For a very general hypersurface $X$ of degree $d$ as considered above, with the field of functions $K$,  the kernel $A_0(X_K)$ of the map $CH_0(X_K)\to\mathbb Z$ is nontrivial.  But then one could ask to determine if for some integer $N$, we have  $NA_0(X_K)=0$. In general, such integer $N$ exists if $X$ is rationally connected (see for example \cite[Prop.11]{CT05}). An upper bound for $N$ was given by Roitman \cite{roitman} (he considered more generally complete intersections in  projective space). Using specialization techniques, Levine and Chatzistimatiou gave a lower bound in the generic case \cite{LCh}.}
\vspace{0.5cm}

\section{Unramified Brauer group and  fibrations in quadrics}\label{fq}

In this section we explain in details how to compute some of the invariants used above, in particular, for fibrations in quadrics, 
from the point of view of unramified cohomology.
For $K$ a field of $char\,K\neq 2$,  one associates the groups of unramified cohomology $H^i_{nr}(K, \mathbb Z/2)$, as some subgroups of the Galois cohomology groups $H^i(K, \mathbb Z/2)$. When $K=k(X)$ is the function field of an algebraic variety $X$ over a field $k$, these groups provide  birational invariants, that can be useful in the study of rationality properties. 
In degree $2$, we essentially obtain the Brauer group (see below). 
In the next sections we will be interested in fibrations in quadrics, when only $\mathbb Z/2$-coefficients 
may lead to nontrivial invariants. 
More generally,  one also considers torsion coefficients $\mu_{n}^{\otimes j}$, for $n$ invertible in $K$, where 
$\mu_{n}$ is the \'etale $k$-group scheme of the $n^{th}$ roots of unity. For a positive integer   $j$  we write 
$$
\mu_{n}^{\otimes j}=\mu_{n}\otimes\ldots\otimes\mu_{n}
$$ 
($j$ times).  If $j<0$, we set 
$$
\mu_{n}^{\otimes j}=Hom_{k-gr}(\mu_{n}^{\otimes (-j)}, \mathbb Z/n)
$$ 
and $\mu_{n}^{\otimes 0}=\mathbb Z/n$.  If $K$  contains a primitive $n^{th}$ root of unity, we have an isomorphism   $\mu_{n}^{\otimes j}\stackrel{\sim}{\to}\mathbb{Z}/n$ for any~$j$.

\subsection{Residues}

Let $A$ be a discrete valuation ring, with field of fractions $K$ and residue field $\kappa (v)$.   For any integers $i>0$, $j$ and for $n$ invertible on $A$ the residue maps  
$$\partial_v^i: H^i(K, \mu_{n}^{\otimes j})\to H^{i-1}(\kappa(v), \mu_{n}^{\otimes (j-1)}) $$  are well-defined. 
Sometimes we will omit $i$ and write just  $\partial_v$.
For $n=2$ we obtain 
$$
\partial_v^i: H^i(K,\mathbb Z/2)\to  H^{i-1}(\kappa(v), \mathbb Z/2).
$$

If $A_v$ is the completion of  $A$ and $K_v$ the field of fractions of $A_v$,  then the residue map factorizes through $K_v$ as follows:
\begin{equation}\label{resfact}\partial_v^i: H^i(K, \mu_{n}^{\otimes j})\to H^i(K_v, \mu_{n}^{\otimes j})\to H^{i-1}(\kappa(v), \mu_{n}^{\otimes (j-1)}). \end{equation}

See \cite{CT, CTOj} for more details, here we just need the cases $i=1,2$ explained below:  

 \begin{enumerate}
 \item $n=2, i=1$. 
 Recall that Kummer theory provides an isomorphism  $$H^1(K, \mathbb Z/2)\simeq K^*/K^{*2}.$$  If $a\in K^*$, we will still denote by $a$ its class in $H^1(K, \mathbb Z/2)$, where no confusion is possible.   We also have $H^0(K, \mathbb Z/2)=\mathbb Z/2$ and $\partial^1_v(a)=v(a)$~mod~$2$. 
\item    From  Kummer theory again, we have $Br(K)[2]=H^2(K,\mathbb Z/2)$. Using cup-products, for $a,b\in K^*$,  we obtain the symbol 
$$
(a,b):=a\cup b\in H^2(K, \mathbb Z/2).
$$  
Then
\begin{equation}\label{resdeg2}\partial^2_v(a, b)=(-1)^{v(a)v(b)}\overline{\frac{a^{v(b)}}{b^{v(a)}}},\end{equation} where $\overline{\frac{a^{v(b)}}{b^{v(a)}}}$ is the image of the unit $\frac{a^{v(b)}}{b^{v(a)}}$ in $\kappa(v)^*/\kappa(v)^{*2}$.\\
\end{enumerate}

\rem{In general, one can explicitly compute the residues on symbols $a_1\cup\dots \cup a_i \in  H^i(K, \mu_{n}^{\otimes i})$, where $a_1,\ldots a_i\in H^1(K, \mu_n)\simeq K^*/K^{*n}$, using the valuations and reductions in the residue field (see  for example \cite{CTOj}). Recall also that any element in   $H^i(K, \mu_{n}^{\otimes i})$ is a sum of symbols: this follows from the Bloch-Kato conjecture (theorem of Voevodsky).\\}

The Gersten conjecture for discrete valuation rings (a theorem, see  \cite[(3.10)]{CT})  describes the elements with trivial residues: in fact, there is an exact sequence:
\begin{equation}\label{Gdvr}
0\to  H^i_{\acute{e}t}(A, \mu_{n}^{\otimes j})\to H^i(K, \mu_{n}^{\otimes j})\stackrel{\partial_v^i}\to H^{i-1}(\kappa(v), \mu_{n}^{\otimes (j-1)})\to 0
\end{equation}

When $B$ is a regular local ring one can look at discrete valuations on $B$ corresponding to height one prime ideals. The following statement  follows from \cite[Thm. 3.82]{CT}:
\propd\label{purity0}{\cite[Thm. 3.82]{CT} \it Let $B$ be a local ring of a smooth variety over a field $k$ and $K$ the field of fractions of $B$. Let $\alpha\in  H^i(K, \mu_{n}^{\otimes j})$, $(n,char\,k)=1$ be such that $\partial^i_v(\alpha)=0$ for all valuations $v$ corresponding to height one prime ideals of $B$. Then $\alpha$ is  in the image of the natural map   $ H^i_{\acute{e}t}(B, \mu_{n}^{\otimes j})\to H^i(K, \mu_{n}^{\otimes j})$.\\}

Note that the map $ H^i_{\acute{e}t}(B, \mu_{n}^{\otimes j})\to H^i(K, \mu_{n}^{\otimes j})$ is in fact injective \cite[Thm. 3.8.1]{CT}. The following corollary will be useful:

\coro\label{purity}{Let $B\subset A$ be an inclusion of local rings, with fields of fractions $K\subset L$ respectively.  Assume that $B$ is a local ring of a smooth variety over a field and that $A$ is a discrete valuation ring with valuation $v$. Let   $\alpha\in  H^i(K, \mu_{n}^{\otimes j})$, $(n,char\, K)=1$ be such that $\partial^i_w(\alpha)=0$ for all valuations $w$ corresponding to height one prime ideals of $B$. Then the image $\alpha'$ of $\alpha$ in  $H^i(L, \mu_{n}^{\otimes j})$  satisfies $\partial^i_v(\alpha')=0$. }
\proof{By Proposition \ref{purity0}, $\alpha$ is   in the image of   $ H^i_{\acute{e}t}(B, \mu_{n}^{\otimes j})\to H^i(K, \mu_{n}^{\otimes j})$, hence $\alpha'$ is in the image of the map    $ H^i_{\acute{e}t}(A, \mu_{n}^{\otimes j})\to H^i(L, \mu_{n}^{\otimes j})$ and the statement follows from  (\ref{Gdvr}). \qed\\}

We also have the following compatibility statement (see \cite[Section 1]{CTOj}):
\propd\label{rescomp}{\it Let $A\subset B$ be discrete valuation rings, with fields of fractions $K\subset L$ respectively. Let $\pi_A$ (resp. $\pi_B$) be the uniformizing parameter of $A$ (resp. $B$). Let $\kappa(A)$ (resp. $\kappa(B)$) be the residue field of $A$ (resp. of $B$).  Let $e$ be the valuation of $\pi_A$ in $B$. We have the following commutative diagram:
$$\xymatrix{
H^i(L, \mathbb Z/2)\ar[r]^{\partial_{B}^i} &H^{i-1}(\kappa(B), \mathbb Z/2)&\\
H^i(K, \mathbb Z/2)\ar[r]^{\partial_{A}^i}\ar[u]^{Res_{K/L}} &H^{i-1}(\kappa(A), \mathbb Z/2)\ar[u]_{eRes_{\kappa(A)/\kappa(B)}} &
}
$$
where $Res_{K/L}$ (resp. $Res_{\kappa(A)/\kappa(B)}$ ) are the restriction maps in Galois cohomology.\\
}

\subsection{Unramified cohomology}

\defi\label{defnr}{For $X$  an integral variety over a field $k$,  integers $j$ and $i\geq 1$ and $n$ invertible in $k$, the {\bf unramified cohomology groups} are defined by
$$\nH^i(X, \mu_n^{\otimes j})=\nH^i(k(X)/k, \mu_n^{\otimes j})=\bigcap\limits_v\mathrm{Ker}[H^i(k(X), \mu_n^{\otimes j})\stackrel{\partial_{v}^i}{\to}H^{i-1}(\kappa(v), \mu_n^{\otimes (j-1)})],$$ where the intersection is over all discrete valuations $v$ on  $k(X)$ (of rank one), trivial on the field $k$.\\}

From the definition, these groups are birational invariants of the variety $X$. 
If $X$ is a smooth projective variety over a field  $k$, then one needs to consider only discrete valuation rings in codimension $1$ points:

$$\nH^i(k(X)/k, \mu_n^{\otimes j})=\bigcap\limits_{x\in X^{(1)}}\mathrm{Ker}[H^i(k(X), \mu_n^{\otimes j})\stackrel{\partial_{x}^i}{\to}H^{i-1}(\kappa(x), \mu_n^{\otimes (j-1)})],$$ where the intersection is over all discrete valuation rings $\mathcal O_{X,x}$ in points $x$ of codimension $1$ of  $X$,  $\kappa(x)$ denotes the residue field of $x$ and  ${\partial_{x}^i}$ is the corresponding residue map.\\

See \cite{CT} for the properties of these groups. We will need the following:
\propd\label{nrst1} {\it 
 If $X$ is a stably rational variety over a field $k$, the natural maps $H^i(k, \mu_n^{\otimes j})\to \nH^i(k(X)/k, \mu_n^{\otimes j})$ are isomorphisms for all $i\geq 1$.\\
}

In degree $2$ we find the Brauer group, so that the point of view of the unramified cohomology could be viewed as a method to exhibit nontrivial elements in this group:

\propd\label{nrst} {\it 
If $X$ is a smooth projective variety over a field $k$, then  $$ \Hn(k(X)/k, \mu_n)\simeq Br(X)[n].$$\\
}
\vspace{-0.5cm}

\rem\label{trhyp}{\begin{enumerate}\item Let $X\subset\mathbb P^n_{\mathbb C}$ be a smooth hypersurface of degree $d\leq n$ with $n\geq 4$. Then $H^3(X, \mathbb Z)_{tors}=0$  and $Br(X)=0$ (see  \cite[Thm. 5.6]{CT15} and \cite[Prop. 4.2.3]{CT}).  Then $ \nH^3(\mathbb C(X)/\mathbb C, \mu_n)=0$ in the following cases (see \cite[Thm. 5.6, Thm. 5.8]{CT15}): $n=4$ or $n>5$;   $n=5$ and $d=3$.
\item  A motivation to study the fibrations in quadrics comes in particular from the fact, that a cubic fourfold containing a plane is birational to such a fibration, with discriminant curve of degree $6$.  By \cite[Thm. 1]{CTAP},  using the point of view of fibrations in quadrics,  for a very general cubic fourfold $X$ over $\mathbb C$, containing a plane, the third unramified cohomology group of $X$ is universally trivial.\\
\end{enumerate}  }

From  Proposition \ref{nrst1}, we have $\Hn (\mathbb C(S)/\mathbb C, \mu_n)=0$ for $S$ a rational complex surface.  We will need a stronger statement (see for example \cite[Thm.1]{artinmumford}):

\propd\label{seqS} {\it Let $S$ be a smooth projective rational complex surface and $K$ 
the function field of $S$. We have the following exact sequence
\begin{equation}\label{exseqS}0\to Br(K)[n]\stackrel{\oplus\partial^2}{\to} \oplus_{x\in S^{(1)}}H^1(\kappa(x), \mathbb Z/n)\stackrel{\oplus\partial^1}{\to}  \oplus_{P\in S^{(2)}} Hom(\mu_n, \mathbb Z/n), 
\end{equation}
where $S^{(r)}$ is the set of points of codimension  $r$ of $S$; $\kappa(x)$  (resp. $\kappa(P)$) is the residue field of $x$ (resp. of $P$),  and the maps are induced by residues. 
}

\vspace{0.7cm}

\subsection{Function fields of quadrics}

Let $K$ be a field, $char\,K\neq 2$, and $q$  a nondegenerate quadratic form over  $K$.  We write
$$q\simeq \langle a_1,\ldots, a_n\rangle$$ for the diagonal form of $q$ in an orthogonal basis, such that
$q(x)=a_1x_1^2+\ldots+a_nx_n^2$. We say that $q$ is similar to a form $ \langle a_1,\ldots, a_n\rangle$ if $q\simeq \langle ca_1,\ldots, ca_n\rangle$ for a constant $c\in K^*$.
 For the general theory of quadratic forms see for example \cite{fqbook}. Here we will be mostly interested in the cases $n=\mathrm{dim}\,q$ is $3$ or $4$. 
We will use the following two cohomological invariants of $q$.
\begin{itemize}
\item The {\bf discriminant} of $q$ is the class  $disc(q)$ of $(-1)^{n(n-1)/2}a_1\ldots a_n$ in $K^*/K^{*2}$. In dimension $4$, any quadratic form $q$ is similar to a quadric of the form 
$$
q\simeq \langle 1,-a,-b,abd\rangle,\quad \text{  where } \quad d=disc(q).
$$ 
\item The {\bf Clifford invariant} of $q$ is an element $c(q)\in H^2(K,\mathbb Z/2)$. In dimension $4$, we have
$$\text{if }q\simeq\langle 1,-a,-b,abd\rangle, \text{ then } c(q)=(a,b)+(ab,d)$$  and $c(\lambda q)=c(q)+(\lambda, disc(q)).$
 \end{itemize}

 Let $Q$ be the  quadric defined by a homogeneous equation $q=0$. The natural maps $H^i(K, \mathbb Z/2)\to  H^i(K(Q), \mathbb Z/2)$ induce the maps 
 
 $$\tau_i:H^i(K, \mathbb Z/2)\to  \nH^i(K(Q)/K, \mathbb Z/2).$$

For $i\leq 4$ these maps have been  studied in \cite{KRS}. In the next section we will need only the case $i\leq 2$, going back to Arason \cite{Arason}.
 
\thmd\label{ffq}{\it \begin{enumerate}
\item If $i=1$ and $q\simeq \langle 1, -a\rangle$, then $ker(\tau_i)\simeq \mathbb Z/2$, generated by the class of $a$. 
\item If $i=1$ and dim\,$q\neq 2$   then the map $\tau_i$ is injective.
\item If $i=2$ and $disc(q)$ is not a square, then $\tau_i$ is a bijection.
\item If $i=2$ and $q\simeq <1,-a,-b,ab>$ is a Pfister form, then   $ker(\tau_i)\simeq \mathbb Z/2$, generated by the class $(a,b)$. 
\end{enumerate}}

\rem{In what follows we often have $\mathbb C\subset K$, so that we omit the minus signs for the coefficients of quadrics.\\}

We will need the following corollary.

\coro\label{resnul}{Let $A$ be a local ring with field of fractions $K$ and residue field $k$  
and $q\simeq <1,-a,-b,abd>$ a quadratic form over $K$. Let  $Q/K$ be the corresponding quadric. Let $v$ be a discrete valuation on $K(Q)$ with $B$ the valuation ring. Assume we have an injection of local rings $A\subset B$. Assume that, up to  multiplication by a square in $K$, the element $d$ is a unit in $A$,  that is a square in $k$. Let $\alpha=(a,b)\in H^2(K, \mathbb Z/2)$ and $\alpha'$ be its image in $H^2(K(Q), \mathbb Z/2)$.  Then $\partial_v(\alpha')=0$.

\proof{Let $\hat A$  be the completion of $A$ and $\hat K$ its field of fractions.  
By hypothesis, over the field $\hat K$ the element $d$ is a square and hence $q\simeq <1,-a,-b,ab>$ over $\hat K$. In particular, $\alpha$ is in the kernel of the map   $H^2(\hat K, \mathbb Z/2)\to H^2(\hat K(Q), \mathbb Z/2)$ by Theorem \ref{ffq}.  But,  using (\ref{resfact}),  the composite map  
$$H^2(K, \mathbb Z/2) \to H^2(\hat K, \mathbb Z/2) \to H^2(\hat K(Q), \mathbb Z/2)\to H^2(K(Q)_v, \mathbb Z/2)\stackrel{\partial}{\to} H^{1}(\kappa(v), \mathbb Z/2)$$
sends $\alpha$ to the residue $\partial_v(\alpha')$, hence  $\partial_v(\alpha')=0$.    \qed}

\vspace{0.7cm}

\subsection{Birational strategy.}
Artin and Mumford \cite{artinmumford}  constructed examples of fibrations in conics  $X$ over $\mathbb P^2_{\mathbb C}$ with $Br(X)[2]\neq 0$. In the work of Colliot-Th\'el\`ene and Ojanguren \cite{CTOj}, these examples are explained using the point of view of unramified cohomology.   
Let us sketch this general birational strategy:
\begin{enumerate}
\item  Let $S$ be a smooth, projective, rational surface over $\mathbb C$, let $K$ be the field of functions of $S$. Let $X$ be a smooth projective variety with a morphism $X\to S$ with generic fiber a smooth conic $Q/K$.  The goal is to compute the group $Br(X)=Br(X)[2]=\Hn(K(Q)/\mathbb C, \mathbb Z/2 )$  or to understand if this group is nonzero.   
 The Hochschild-Serre spectral sequence implies that the natural map $Br(K)\to Br(Q)$ is surjective.  
One then uses that $Br(Q)=Br(Q)[2]=\Hn(K(Q)/K, \mathbb Z/2 )$ and  the straightforward  inclusion  
 $ \Hn(K( Q)/K, \mathbb Z/2 )\supset \Hn(K(Q)/\mathbb C, \mathbb Z/2 )$
to claim that any element $\xi$ in the  latter group comes from an element $\beta\in Br(K).$
\item  The description of $\beta$ via family of its residues, is given by the exact sequence (\ref{exseqS}). To show that an element  $\beta$ gives an element in 
$Br(X)=\Hn(K(Q)/\mathbb C, \mathbb Z/2 )$, we need to show that for any  discrete valuation $v$ on $K(Q)$ the residue  $\partial_v(\beta)$ is zero. Let $A$ be the valuation ring of $v$. 
The rational map $\mathrm{Spec}\,A\dashrightarrow K$ induced by the inclusion $K\subset K(Q)$ extends to a morphism 
$\mathrm{Spec}\,A \to S$ and the image $x_v$ of the closed $\kappa(v)$-point of 
$\mathrm{Spec}\,A$ is called the {\bf center} of $v$ in $S$. There are two cases to consider: $x_v$ is the generic point of  curve $C_v$  or a closed point $P_v$. This is the most technical part, done case-by-case using Proposition \ref{rescomp} and other local arguments.\\
\end{enumerate}

For complex threefolds fibered in conics over a  rational surface, Colliot-Th\'el\`ene (Lecture course, Beijing BICMR, 2015) has obtained a general formula for the Brauer group of $X$:

\thmd{\it (Colliot-Th\'el\`ene)
Let $S$ be a smooth, projective, rational surface over $\mathbb C$ and $K$ the field of functions of $S$. Let $X$ be a smooth threefold equipped with a conic bundle structure $\pi: X\to S$ and $\alpha \in Br(K)[2]$ 
the class corresponding to the quaternion algebra associated to the conic given by the generic fiber of $\pi$. Assume that $\alpha$ is nonzero and that the ramification curve $C$, consisting of codimension $1$ points $x$ of $S$ such that $\partial_x(\alpha)\neq 0$, is with only quadratic singularities. Let $C=\cup_{i=1}^n\, C_i\subset S$ be its decomposition into irreducible components and $(\gamma_i)$ the associated family of residues of $\alpha$ in $\oplus_{i=1}^n H^1(\kappa(C_i),\mathbb Z/2)$. Consider the following subgroup $H\subset (\mathbb Z/2)^n$:
$$H=\{(n_i)\,|\, n_i=n_j \text{ for }i\neq j, \text{  if there is a point }P\in C_i\cap C_j, \partial_P(\gamma_i)=\partial_P(\gamma_j)\neq 0\}.$$
Then the Brauer group $Br(X)$ is the quotient of $H$ by the diagonal $(1,\ldots ,1) \mathbb Z/2$.\\} 

\rem\label{trirr}{In the theorem above, if $C$ is irreducible, then  $Br(X)=0$.\\}

Below we provide a similar formula for fibrations in quadrics of dimension $2$, using the birational strategy above. We start with two concrete examples before giving the general description.\\

\subsection{Examples in dimension $4$}
For the projective plane $\mathbb P^2_{\mathbb C}$, with homogeneous coordinates $[X:Y:Z]$ we use standard coordinates $x$ and $y$ (resp. $y$ and $z$, resp. $x$ and $z$) for the open $Z\neq 0$ (resp. $X\neq 0$, resp. $Y\neq 0$). 
Consider a fibration in quadrics over $\mathbb P^2_{\mathbb C}$ with generic fiber $Q/K$, $K=\mathbb C(x,y)$ defined by a quadratic form 
$$
q=\langle x,y,1, xy((x+y+1)^4+xy)\rangle,
$$ 
i.e.,
\begin{equation}\label{zn0}
Q: xx_0^2+yx_1^2+x_2^2+xy((x+y+1)^4+xy)x_3^2=0.\end{equation}
 Note that   over an open $X\neq 0$ we have a model of $Q$ defined by
\begin{equation}\label{xn0}
x_0^2+yx_1^2+zx_2^2+yz((1+y+z)^4+yz^2)x_3^2=0.\end{equation}

The class of the discriminant of $q$ in $K^*/K^{*2}$ is the class of $d=((x+y+1)^4+xy)$, hence nonzero. 

\propd\label{brntr}{\it Let $\alpha=(x, y)\in Br(K)[2]$. Then the image $\alpha'$ of  $\alpha$ in $H^2(K(Q), \mathbb  Z/2)$ is a nonzero class, that lies in the  unramified cohomology subgroup $$\alpha'\in H^2_{nr}(K(Q)/\mathbb C, \mathbb Z/2)\subset H^2(K(Q), \mathbb  Z/2).$$}

To establish the proposition we have to show that for any discrete valuation $v$ on $K(Q)$ we have $\partial_v(\alpha')=0$. Note that by Theorem \ref{ffq}, the class $\alpha'$ is nonzero.
Let us first investigate the ramification of $\alpha$ on $\mathbb P^2_{\mathbb C}$. From the definition and the formula (\ref{resdeg2}), we  only have the following nontrivial residues:
\begin{itemize}
\item  $\partial_x(\alpha)=y$ at the line $L_x: x=0$, where we write $y$ for its class in the residue field $\mathbb C(y)$ modulo squares;
\item $\partial_y(\alpha)=x$ at the line $L_y: y=0$,
\item $\partial_z(\alpha)=\partial_z(z, zy)=y$ at the line $L_z: z=0$, where we view the generic point of $L_z$ in the open $X\neq 0$, with homogeneous coordinates $y$ and $z$. 
\end{itemize}
Let $A$ be the valuation ring and $x_v$ be the center of $v$ in  $\mathbb P^2_{\mathbb C}$. There are two cases to consider: $x_v$ is the generic point of  a curve $C_v$  or a closed point $P_v$.

\subsubsection{Codimension $1$ case}
Recall that by Proposition \ref{rescomp}, using the notations above, the inclusion of discrete valuation rings 
$\mathcal O_{\mathbb P^2, C_v}\subset A$ induces a commutative diagram
\begin{equation}\label{diagr}\xymatrix{
H^2(K(Q), \mathbb Z/2)\ar[r]^{\partial^2} &H^{1}(\kappa(v), \mathbb Z/2)&\\
H^2(K, \mathbb Z/2)\ar[r]^{\partial^2}\ar[u] &H^{1}(\kappa(C_v), \mathbb Z/2)\ar[u] &
}
\end{equation}
Hence we have the following cases:
\begin{enumerate}
\item $C_v$ is different from  $L_x, L_y$ or $L_z$. Then $\partial_{C_v}(\alpha)=0$, so that $\partial_v(\alpha')$ is zero from the diagram above.
\item $C_v$ is one of the lines  $L_x, L_y$ or $L_z$. Then note that modulo the equation of $C_v$ the element  $d=(x+y+1)^4+xy$ (or $(X+Y+Z)^4+XYZ^2$ in homogeneous coordinates) is a nonzero square, so that Corollary \ref{resnul} gives $\partial_v(\alpha')=0$.
\end{enumerate}
We deduce that for any valuation $v$ on $K(Q)$ with center a codimension $1$ point in $\mathbb P^2_{\mathbb C}$ the residue $\partial_v(\alpha')$ is zero.

\subsubsection{Codimension $2$ case} Let $P_v$ be the center of $v$ on $\mathbb P^{2}_{\mathbb C}$. We have an inclusion of local rings  $\mathcal O_{\mathbb P^2, P_v}\subset A$ inducing the inclusion of corresponding completions 
$\widehat{\mathcal O_{\mathbb P^2, P_v}}\subset  A_v$ with function fields $K_{P_v}\subset K(Q)_v$ respectively. We have three possibilities for the point $P_v$:
\begin{enumerate}
\item If $P_v$ is not on the union $L_x\cup L_y\cup L_z$, then $\alpha$ is a cup product of units in $\mathcal O_{\mathbb P^2, P_v}$, hence units in $A_v$, so that $\partial_v(\alpha')=0$.
\item $P_v$ is on only one curve, say $L_x$. Then the image of $y$ in $\kappa(P_v)$ is a nonzero complex number, hence a square in   $\widehat{\mathcal O_{\mathbb P^2, P_v}}$, hence $y$ is also a square in $A_v$. We deduce that $\alpha'=0$ in $H^2(K(Q)_v, \mathbb Z/2)$, hence $\partial_v(\alpha')=0$.
\item $P_v$ is on two curves, say $P=L_x\cap L_y$. Then the image of $d=(x+y+1)^4+xy$ in $\kappa(P_v)$ is a nonzero complex number, hence a square in   $\kappa(P_v)$.
Then Corollary \ref{resnul} gives $\partial_v(\alpha')=0$.

\end{enumerate}

We deduce that for any valuation $v$ on $K(Q)$ with center a codimension $2$ point in $\mathbb P^2_{\mathbb C}$ the residue $\partial_v(\alpha')$ is zero. This concludes the proof of the proposition.

\vspace{1cm}

Let us consider another example. Let $C_i, i=1,2,3,$ be smooth conics in 
$\mathbb P^2_{\mathbb C}$ given by equations $f_i(X,Y,Z)=0$, such that $C_1$ and $C_2$ intersect in $4$ distinct points $P_1, \ldots, P_4$, the conic $C_3$ is bitangent to $C_1$ at points $P_5$ and $P_6$ (distinct from previous four points) and $C_3$ is bitangent to $C_2$ at points $P_7$ and $P_8$ (distinct from previous six points). For example, one could take
$$f_1=\frac{X^2}{4}+Y^2-Z^2, \, f_2=X^2+\frac{Y^2}{4}-Z^2, \,f_3=X^2+Y^2-Z^2.$$ 

Consider a fibration in quadrics over $\mathbb P^2_{\mathbb C}$ with generic fiber $Q/K$, $K=\mathbb C(x,y)$ defined by a quadratic form $q=\langle 1,f_1,f_2, f_1f_2f_3\rangle$.
The class of the discriminant of $q$ in $K^*/K^{*2}$ is the class of $f_3$, hence nonzero. 
\propd\label{brntr2}{\it Let $\alpha=(f_1, f_2)\in Br(K)[2]$. Then the image $\alpha'$ of  $\alpha$ in $H^2(K(Q), \mathbb  Z/2)$ is a nonzero class, that lies in the  unramified cohomology subgroup $$\alpha'\in H^2_{nr}(K(Q)/\mathbb C, \mathbb Z/2)\subset H^2(K(Q), \mathbb  Z/2).$$}

As before, let $v$ be a valuation on $K(Q)$ with $A$  the valuation ring and $x_v$  the center of $v$. 
\subsubsection{Codimension $1$ case} Assume  $x_v$ is the generic point of  a curve $C_v$. 
\begin{enumerate}
\item $C_v$ is different from  $C_i$, $i=1,2$. Then $\partial_{C_v}(\alpha)=0$, so that $\partial_v(\alpha')$ is zero from the diagram of proposition \ref{rescomp}.
\item $C_v$ is one of the two conics  $C_1, C_2$. Then the bitangency   condition implies that, modulo the equation  of $C_v$, the element  $f_3$ is a nonzero square, so that Corollary \ref{resnul} gives $\partial_v(\alpha')=0$.
\end{enumerate}

\subsubsection{Codimension $2$ case} Let $P_v$ be the center of $v$ on $\mathbb P^{2}_{\mathbb C}$. We have an inclusion of local rings  $\mathcal O_{\mathbb P^2, P_v}\subset A$ inducing the inclusion of corresponding completions 
$\widehat{\mathcal O_{\mathbb P^2, P_v}}\subset  A_v$. 
\begin{enumerate}
\item If $P_v$ is not on the union $C_1\cup C_2$, then $\alpha$ is a cup product of units in $\mathcal O_{\mathbb P^2, P_v}$, hence units in $A_v$, so that $\partial_v(\alpha')=0$.
\item $P_v$ is on only one curve, say $C_1$. Then the image of $f_2$ in $\kappa(P_v)$ is a nonzero complex number, hence a square in   $\widehat{\mathcal O_{\mathbb P^2, P_v}}$, hence $f_2$ is also a square in $A_v$. We deduce that $\alpha'=0$ in $H^2(K(Q)_v, \mathbb Z/2)$, hence $\partial_v(\alpha')=0$.
\item $P_v=P_i$, $i=1,\ldots , 4$. Then the image of $f_3$ in $\kappa(P_v)$ is a nonzero complex number, hence a square in   $\kappa(P_v)$.
Then corollary \ref{resnul} gives $\partial_v(\alpha')=0$.

\end{enumerate}

We deduce that for any valuation $v$ on $K(Q)$ the residue $\partial_v(\alpha')$ is zero, that finishes the proof of proposition.

Note that in this example the divisor $C_1\cup C_2\cup C_3$ is not in simple normal crossings, but the construction still works.

\vspace{0.5cm}

\subsection{Brauer group of fibrations in quadric surfaces}
Let  $S$ be a smooth, projective, rational surface over $\mathbb C$, let $K$ be the field of functions of $S$. Let $X$ be a fourfold equipped with a projective map $\pi:X\to S$ with generic fiber a quadric $Q/K$ given by a nondegenerate quadratic form $q$ over $K$. First observe:
\begin{itemize}
\item The quadric $Q$ becomes rational over (at most) a degree two extension of $K$, 
hence, by a restriction-corestriction argument, the Galois cohomology of $K(Q)$ is $2$-torsion, 
so that we are only interested in the group 
$\Hn(K(Q)/\mathbb C,\mathbb Z/2)$, isomorphic to the Brauer group of $X$, if $X$ is smooth.
\item If the discriminant $d=disc(q)$ is a square, then the quadric $Q$ is $K$-birational to a product  $C\times C$ for $C$ a conic over $K$ \cite[Thm. 2.5]{CTSk}. More precisely, if $q$ is similar to a form $\langle 1,-a,-b,ab\rangle$ then one could take for $C$ the conic defined by  the form $\langle 1,-a,-b\rangle$. Let $L=K(C)$. Then the field $K(Q)=L(C_L)$ is a purely transcendental extension of $L$ since $C(L)\neq \emptyset$.  Hence this case is essentially equivalent to the conic bundles case, so that we will only consider the case when $disc(q)$ is nontrivial (i.e. nonsquare).
\end{itemize}

In the  statement below
we use the following notation: for $C\subset S$ an irreducible curve we write $\partial_C$ for the residue at the generic point of $C$;  we also identify the divisor $D=\sum D_i$ with the set of generic points of $D_i$.

\thmd\label{BrQBS}{\it Let  $S$ be a smooth, projective, rational surface over $\mathbb C$ and $K$ the field of functions of $S$. Let  $Q/K$ be a two-dimensional quadric given by a nondegenerate quadratic form $q$ over $K$. 
Let $\alpha\in Br(K)$ be the Clifford invariant $\alpha=c(q)$ of $q$. Assume that the divisor 
$$
ram\, \alpha=\{x\in S^{(1)}\,|\,  \partial_x(\alpha)\neq 0\}
$$ 
is a simple normal crossing divisor. Let $\mathcal P$ be the set of singular points  of $ram\, \alpha$.  Let $d\in K^*/K^{*2}$ be the discriminant of $q$. Assume that $d$ is nontrivial and consider the following divisor:
\begin{multline*} T=\{x\in S^{(1)}\,|\, \partial_x(\alpha)\neq 0 \text{ and }  d, \text{ up to a multiplication by a square}, \\ \text{is a unit in } \mathcal O_{S,x}, \text{ and the image of }d\text{ in }\kappa(x)\text{ is a square.}\}
\end{multline*}
Let $T=\cup_{i=1}^n T_i$ be the decomposition into irreducible components. For $i=1,\ldots, n$ 
put $c_i:=\partial_{T_i} (\alpha)$.
Let 
\begin{equation}\label{defH}
H=ker[(\mathbb Z/2)^n \stackrel{\partial^1}{\to} \oplus_{P\in \mathcal P} H^0(\kappa(P), \mathbb Z/2)],\; (n_i)_{i=1}^n\mapsto (\oplus n_i\,\partial^1_P (c_i))
\end{equation}
Then  the natural homomorphism
$$H\to H^2(K(Q),\mathbb Z/2),$$
that associates to a family $(n_i)_i\in H$ the image $\beta'$ in  $H^2(K(Q),\mathbb Z/2)$ of the unique class $\beta\in H^2(K,\mathbb Z/2)$ with 
$$
\partial^2(\beta)=(n_ic_i)_i\in \oplus_{i=1}^n H^1(\kappa(T_i),\mathbb Z/2),
$$ 
via the exact sequence (\ref{exseqS})
$$0\to Br(K)[2]\stackrel{\partial^2}{\to} \oplus_{x\in S^{(1)}}H^1(\kappa(x), \mathbb Z/2)\stackrel{\partial^1}{\to}  \oplus_{P\in S^{(2)}} H^0(\kappa(P), \mathbb Z/2)$$ 
induces an isomorphism
$$\Phi: H\stackrel{\sim}{\to}  \Hn(K(Q)/\mathbb C,\mathbb Z/2).$$

}

\vspace{0.5 cm}

\rem\label{remexplicite}{\begin{enumerate}
\item In the construction $\partial_x(\beta)=0$ if $x$ is not the generic point of one of the components of $T$.
\item For a given quadric surface $Q/K$ we could always assume that $ram\, \alpha$ is a simple normal crossings divisor, 
after blowing up $S$.
\item Explicitly, up to similarity,  we may assume that $q$ can be represented as $q=\langle 1, a,b, abd\rangle$ (note that we can ignore the signs since we are over $\mathbb C$), so that $c(q)=(a,b)+(ab,d)$. If $x\in T$, then 
\begin{equation}\label{dcq}
\partial_x(\alpha)=\partial_x(a,b)
\end{equation}
\item The definition of $H$ as the kernel of the map (\ref{defH}) is precisely to ensure that $H$ is in the kernel of the map $\partial^1$ of the sequence (\ref{exseqS}).  
\end{enumerate}}

\vspace{0.5cm}

The proof of Theorem \ref{BrQBS} follows the same strategy as in the preceding examples. We split into two parts: first we show that any element of $H$ gives an unramified element on $K(Q)$ and then we establish that any element in $\Hn(K(Q)/\mathbb C,\mathbb Z/2)$ comes from $H$. Note that $H\to H^2(K(Q),\mathbb Z/2)$
 is an injection by Theorem \ref{ffq}.
 
\subsubsection{$\Phi$ is well-defined.}
Fix $(n_i)_i\in H$ and let $\beta\in H^2(K,\mathbb Z/2)$ be such that 
$$
\partial^2(\beta)=(n_ic_i)_i\in \oplus_{i=1}^n H^1(\kappa(T_i),\mathbb Z/2)
$$ 
and $\beta'$ be its image in $H^2(K(Q),\mathbb Z/2)$. Note that the conditions on $n_i$'s guarantee that $(n_ic_i)_i$ is in the kernel of $\partial^1$ in the sequence (\ref{exseqS}),  so that $\beta$ exists.  We have to show that for any discrete valuation $v$ on $K(Q)$ we have $\partial_v(\beta')=0$.
Let $A$ be the valuation ring and $x_v$ be the center of $v$ in  $S$. There are two cases to consider: $x_v$ is the generic point of  a curve $C_v$  or a closed point $P$.

{\it Codimension $1$ case.} If   $\partial_{C_v}(\beta)=0$, then  $\partial_v(\beta')$ is zero by Proposition \ref{rescomp}. If  $\partial_{C_v}(\beta)\neq 0$, then $C_v$ is one of the curves $T_i$ and by definition  $\partial_{C_v}(\beta-\alpha)=0$. But, because of the condition on $d$ we  could apply Corollary \ref{resnul} to get  $\partial_v(\alpha')=0$, where $\alpha'$ is the image of $\alpha$ in $H^2(K(Q), \mathbb Z/2)$. Hence $\partial_v(\beta')=0$ as well. \\
{\it Codimension $2$ case.}
\begin{enumerate}
\item If $P$ is not on $T$, then $\beta$ is unramified for the valuations corresponding to the height one prime ideals of $\mathcal O_{S,P}$, so that $\partial_v(\beta')=0$ by Corollary \ref{purity}.
\item $P\in T_i\cap T_j$ with $n_i=n_j=1$ and  $\partial_P(c_i)=\partial_P(c_j)\neq 0$.  Then $\beta-\alpha$ is unramified for the valuations corresponding to the height one prime ideals of  $\mathcal O_{S,P}$ and  $\partial_v(\beta'-\alpha')=0$  by Corollary \ref{purity}.
 The image of $d$ in $\kappa(P)$ is a nonzero complex number, hence a square in   $\kappa(P)$.
Then Corollary \ref{resnul} gives $\partial_v(\alpha')=0$, so that $\partial_v(\beta')=0$ as well.
\item In the remaining cases, let $P\in T_i$ with $n_i\neq 0$ (if $P$ is on another curve $T_j$ then $\partial_P(c_j)=0$ by definition of $H$). Then   $c_i=\partial_{T_i}(\beta')$ could be written (modulo squares) as $c_i=u$ with $u$ a unit in $\mathcal O_{T_i, P}$ (using the sequence (\ref{exseqS}) and (\ref{Gdvr})).  Let $w\in \mathcal O_{S,P}$ be such that the image of $w$ in $\kappa(T_i)$ is $u$, then $w$ is nonzero in $\kappa(P)$ and hence a square in $\widehat{\mathcal O_{S,P}}$.
 Let  $t$ be a  uniformizing parameter  of $T_i$ in $K$.
 By construction, $\beta-(t, w)$ is unramified on $\mathcal O_{S,P}$, so that $\partial_v(\beta'-(t, w))=0$ by Corollary \ref{purity}. Since $w$ is a square in $\widehat{\mathcal O_{S,P}}$, and hence in $A$, we deduce $\partial_v(\beta')=\partial_v(t, w)=0$.
\end{enumerate}

\subsubsection{$\Phi$ is surjective.}
By Theorem \ref{ffq} the map  $ H^2(K,\mathbb Z/2)\to \Hn (K(Q)/K,\mathbb Z/2)$ is an isomorphism. Hence any element from $\Hn (K(Q)/\mathbb C,\mathbb Z/2)\subset \Hn (K(Q)/K,\mathbb Z/2)$ comes from an element  $\beta\in  H^2(K,\mathbb Z/2)$. We only need to show that for any $x\in S^{(1)}$ we have 
\begin{equation}\label{check}\partial_x(\beta)\neq 0\Leftrightarrow x=T_i\text{ and }\partial_x(\beta)=c_i.
\end{equation}
The condition at points $P\in T_i\cap T_j$ is automatically forced by the exact sequence (\ref{exseqS}). As in remark \ref{remexplicite}, we choose $q\simeq \langle 1, a,b, abd\rangle$, so that $c(q)=(a,b)+(ab,d)$ and  
$c_i=\partial_{T_i}(a,b).$

So consider $x\in S^{(1)}$ and fix a uniformizing parameter $\pi$ of $A=\mathcal O_{S,x}$, we denote by $v_x$ the corresponding valuation on $K$. Up to changing $q$ by $\lambda q$ with $\lambda\in K^*$ (that does not change $K(Q)$), we may assume that $q$ is with coefficients in $A$ and  $N=0, 1$ or $2$ coefficients of $q$   have odd valuation at $x$. We then have that $q=0$ defines a closed subscheme $Z\subset \mathbb P^4_A$. The divisor $\pi=0$ of $Z$ is  either integral, or a union of two planes. In the former case it defines a discrete valuation $v$ on $K(Q)$, in the latter case we take $v$ corresponding to any of the planes. In Proposition \ref{rescomp} the factor $e=1$ and we have the  following commutative diagram:
$$\xymatrix{
H^2(K(Q), \mathbb Z/2)\ar[r] &H^{1}(\kappa(v), \mathbb Z/2)&\\
H^2(K, \mathbb Z/2)\ar[r]\ar[u]^{Res_{K/K(Q)}} &H^{1}(\kappa(x), \mathbb Z/2)\ar[u]_{Res_{\kappa(x)/\kappa(v)}} &
}
$$
If $N=0$ or $1$, then $\kappa(v)$ is the function field of a quadric  over $\kappa(x)$, that we denote $\bar q=0,$ of positive dimension, so that the map $Res_{\kappa(x)/\kappa(v)}$ is injective by theorem \ref{ffq}. Since the diagram commutes, $\partial_x(\beta)\in ker Res_{\kappa(x)/\kappa(v)}$ and hence trivial for $N=0$ or $1$.
Similarly, if $N=2$, then $\partial_x(\beta)=0$ or $disc(\bar q)$. So the only possible nontrivial residues of $\beta$ correspond to the following cases (up to multiplication by squares):
\begin{itemize}
\item $v_x(a)=1$, $v_x(b)=0$ and $v_x(d)=0$ (or the symmetric case  $v_x(a)=0$, $v_x(b)=1$ and $v_x(d)=0$). Write $a=\pi a_1$ with $a_1$ a unit and $\bar a_1$ (resp. $\bar b$, $\bar d$) for the image in $\kappa(x)$. Then  we have two possibilities for $\bar q$: $\langle 1, \bar b\rangle$ or $\langle \bar a_1, \bar a_1\bar b\bar d \rangle$. If $\partial_x(\beta)\neq 0$ then from the commutative diagram above we must have $\bar b=\bar b\bar d$ and $\partial_x(\beta)=\bar b$, so that $\bar d$ is a square and  $\partial_x(\beta)=\partial_x(a\cup b)$.
\item  $v_x(a)=1$, $v_x(b)=1$ and $v_x(d)=0$. Write $a=\pi a_1$  and $b=\pi b_1$ with  $a_1, b_1$ units as before. Then  we have two possibilities for $\bar q$: $\langle 1, \bar a_1\bar b_1\bar d\rangle$ or $\langle \bar a_1, \bar b_1\rangle$. If $\partial_x(\beta)\neq 0$, then from the commutative diagram above we must have $\bar a_1\bar b_1\bar d=\bar a_1\bar b_1$ and $\partial_x(\beta)=\bar a_1\bar b_1$, so that $\bar d$ is a square and  $\partial_x(\beta)=\partial_x(a\cup b)$.
\end{itemize}
We deduce that the condition (\ref{check}) is satisfied, so that $\beta\in H$. 
This finishes the proof of Theorem \ref{BrQBS}.\qed\\

\coro{Let $\pi: X\to S$ be a fibration in quadrics over a smooth projective rational surface $S$ over $\mathbb C$, such that the total space $X$ is smooth and projective, of dimension $4$. If the discriminant curve $D\subset S$  is smooth and connected, then $Br(X)=0$. 
}
\proof{Let 
$$
\xi\in Br(X)=Br(X)[2]=\Hn(\mathbb C(X)/\mathbb C,\mathbb Z/2)=\Hn(K(Q)/\mathbb C,\mathbb Z/2),
$$ 
where $K$ is the field of functions of $S$ and $Q$ is the generic fiber. As in the proof of the surjectivity above, $\xi$ comes from an element $\beta\in H^2(K, \mathbb Z/2)$ and $\beta$ can have a nontrivial residue only at codimenion one points $s\in S^{(1)}$, such that the fiber $X_s$ is not smooth, since only at the generic point of $D$. However,  at this point, the class of the discriminant of $Q$ has valuation $1$, and hence the set $H$ in the theorem above is empty. \qed\\}

\rem{A cubic fourfold containing a plane is birational to a fibration in quadrics $\pi: X\to S$ with discriminant curve $D$ of degree $6$.  Assume that in 
the construction in Theorem \ref{BrQBS}  the discriminant curve $D$ is of degree $6$ and the class of the discriminant of the generic fiber of $\pi$ is  not a square, then $D=D_0+2D_1$, where $D_0$ corresponds to the discriminant. In particular, $D_1$ is of degree $1$ or $2$. If we had a nontrivial $\beta\in \Hn(\mathbb C(X)/\mathbb C,\mathbb Z/2)$, then, by Theorem \ref{BrQBS},  $\beta$ could have nontrivial residues only along components of $D_1$, via sequence (\ref{exseqS}).  Since the degree of $D_1$ is $1$ or  $2$, one easily checks that this is impossible. }

$\,$

\vspace{1cm}


\vspace{1cm}

\begin{thebibliography}{60}
\bibitem{Arason} J. Kr. Arason, {\it Cohomologische Invarianten quadratischer Formen}, J. Algebra {\bf 36} (1975), p.448--491.

\bibitem{artinmumford} M. Artin et D. Mumford, {\it Some elementary examples of unirational varieties
which are not rational}, Proc. London Math. Soc. (3) {\bf 25} (1972) 75--95.

\bibitem{CTAP} A. Auel, J.-L. Colliot-Th\'el\`ene, R. Parimala, {\it Universal unramified cohomology of cubic
fourfolds containing a plane},    \url{http://arxiv.org/abs/1310.6705v2}, to appear in " Brauer groups and obstruction problems: moduli spaces and arithmetic (Palo Alto, 2013)", Asher Auel, Brendan Hassett, Tony V\'arilly-Alvarado, and Bianca Viray, eds.


\bibitem{survolurnr} A. Beauville, {\it Vari\'et\'es rationnelles et unirationnelles}, Algebraic Geometry -- Open problems (Proc. Ravello 1982), Lecture Notes {\bf 997}, 16--33; Springer-Verlag (1983). 

\bibitem{B6} A. Beauville, \emph{A very general sextic double solid is not stably rational},  arXiv:1411.7484, to appear in Bull. London Math. Soc. 

\bibitem{BeauvilleRapport} A. Beauville, {\it The L\"uroth problem}, arXiv:1507.02476.

\bibitem{BCTSSD} A. Beauville,  J.-L. Colliot-Th\'el\`ene, J.-J. Sansuc, Sir Peter Swinnerton-Dyer, \emph{Vari\'et\'es stablement rationnelles non rationnelles}. Ann. of Math, {\bf 121} (1985) 283--318.



\bibitem{BlochSrinivas} S. Bloch,  V. Srinivas, {\it  Remarks on correspondences and algebraic cycles}. Am. J. Math.
{\bf 105} (1983),  1235--1253.

\bibitem{Bogomolov} F. A. Bogomolov, {\it The Brauer group of quotient spaces of linear representations}, Izv. Akad. Nauk SSSR Ser. Mat. {\bf 51} (1987), no. 3, 485--516, 688;  English transl., Math. USSR-Izv.
{\bf 30}
(1988), no. 3, 455--485.

\bibitem{LCh} M. Levine, A. Chatzistimatiou, {\it Torsion indices of smooth projective varieties}, in preparation.

\bibitem{Cheltsov} I. Cheltsov,  {\it Birationally rigid Fano varieties}, Russian Math. Surveys {\bf 60} (2005), no. 5, 71--160.  

\bibitem{CG72} H. Clemens and P. Griffiths, {\it The intermediate Jacobian of the cubic threefold}, Ann. of Math {\bf 95} (1972) 281--356.

\bibitem{CT} J.-L. Colliot-Th\'el\`ene, \emph{Birational invariants, purity and the Gersten conjecture},  $K$-theory and algebraic geometry: connections with quadratic forms and division algebras (Santa Barbara, CA, 1992),  1--64,
Proc. Sympos. Pure Math., \textbf{58}, Part 1, Amer. Math. Soc., Providence, RI, 1995.

\bibitem{CT05} J.-L. Colliot-Th\'el\`ene, \emph{Un th\'eor\`eme de finitude pour le groupe de Chow des z\'ero-cycles d'un groupe alg\'ebrique lin\'eaire  sur un corps $p$-adique},  Inventiones mathematicae {\bf 159} (2005) 589--606.  

\bibitem{CT15} J.-L. Colliot-Th\'el\`ene, \emph{Descente galoisienne sur le second groupe de Chow : mise au point et applications},  Documenta Mathematica, Extra Volume: Alexander S. Merkurjev's Sixtieth Birthday (2015) 195--220.


\bibitem{CTOj} J.-L. Colliot-Th\'el\`ene, M. Ojanguren,  {\it Vari\'et\'es unirationnelles non rationnelles : au-del\`a de l'exemple d'Artin et Mumford, }Invent. math. {\bf 97} (1989) 141--158.

\bibitem{CTP15} J.-L. Colliot-Th\'el\`ene,  A. Pirutka, \emph{Hypersurfaces quartiques de dimension 3 : non rationalit\'e stable},   Ann. Sci. \'Ecole Norm. Sup. (2) {\bf 49}  (2016) 373--399.

\bibitem{CTP15-2} J.-L. Colliot-Th\'el\`ene,  A. Pirutka, {  \cyrit{Ciklicheskie nakrytiya, kotorye ne yavlyayut\-sya 
  stabil\cprime no  racional\cprime nymi}},  arXiv:1506.00420,   to appear in Izvestiya RAN, Ser. Math.
  


\bibitem{CTSan79} J.-L. Colliot-Th\'el\`ene, J.-J. Sansuc, \emph{La descente sur les vari\'et\'es rationnelles},  Journ\'ees de g\'eom\'etrie alg\'ebrique d'Angers (Juillet 1979), \'edit\'e par A.Beauville, Sijthof et Noordhof (1980) pp. 223--237.

\bibitem{CTSa07} J.-L. Colliot-Th\'el\`ene, J.-J. Sansuc, {\it The rationality problem for fields of invariants under linear algebraic groups (with special regards to the Brauer group)}, Algebraic groups and homogeneous spaces, 113--186, Tata Inst. Fund. Res. Stud. Math., {\bf 19}, Tata Inst. Fund. Res., Mumbai, 2007.

\bibitem{CTSk}  J.-L. Colliot-Th\'el\`ene, A.N. Skorobogatov, {\it Groupes de Chow des z\'ero-cycles des fibr\'es en quadriques}, Journal of K-theory {\bf 7} (1993) 477--500.



  
\bibitem{EZ}  F. El Zein, {\it La classe fondamentale d'un cycle}, Comp. Math. {\bf 29} (1974) 9--33.

\bibitem{dF} T. de Fernex, {\it Birationally rigid hypersurfaces}, Invent. Math. {\bf 192} (2013), 533--566,
Erratum, Invent. Math. {\bf 203} (2016), 675--680.

\bibitem{fulton} W. Fulton, Intersection theory, Springer-Verlag, Berlin, 1998.

\bibitem{HKT15} B. Hassett, A. Kresch, Y. Tschinkel, \emph{ Stable rationality and conic bundles}, arXiv:1503.08497, to appear in Math. Ann.

\bibitem{HT16} B. Hassett,  Y. Tschinkel, \emph{On stable rationality of Fano threefolds and del Pezzo fibrations},  	arXiv:1601.07074. 

\bibitem{IM} V. A. Iskovskikh et Yu. I. Manin, {\it Three-dimensional quartics and counterexamples to the L\"uroth problem}, (in Russian) Mat. Sb. (N.S.) {\bf 86} (128) (1971), 140--166.

\bibitem{KRS} B. Kahn, M. Rost,   R. Sujatha, {\it Unramified cohomology of quadrics, I},
Amer. J. Math. {\bf 120} (1998), 841--891.

\bibitem{fqbook} B. Kahn, {\it Formes quadratiques sur un corps},
Cours sp\'ecialis\'e no {\bf 15}, SMF, 2008. 

\bibitem{kollar} J. Koll\'ar, \emph{Rational curves on algebraic varieties}, Springer-Verlag, Berlin, 1996.

\bibitem{K00}  J. Koll\'ar, \emph{Unirationality of cubic hypersurfaces}, J. Inst.  Math. Jussieu {\bf 1} (2002), 467--476.


\bibitem{Le} F. Lecomte, \emph{Rigidit\'e des groupes de Chow}, Duke Math. J. {\bf 53} (1986) 405-426.


\bibitem{madore} D. Madore, {\it Sur la sp\'ecialisation de la R-\'equivalence}, 
{ \url{http://perso.telecom-paristech.fr/~madore/specialz.pdf}}

\bibitem{Manin}  Yu. I. Manin,
\emph{Cubic forms: algebra, geometry, arithmetic}, Izdat. "Nauka",
Moscow, 1972.

\bibitem{Merkurjev1} A.  S.  Merkurjev, {\it Unramified cohomology of classifying varieties for classical simply
connected groups}, Ann. Sci. \'Ecole Norm. Sup. (4)
{\bf 35}
(2002), no. 3, 445--476.

\bibitem{Merkurjev2} A.  S.  Merkurjev, {\it Unramified elements in cycle modules}, J. London Math. Soc. (2) {\bf 78} (2008) 51--64.

\bibitem{Merkurjev}  A. S. Merkurjev, {\it Invariants of algebraic groups and retract rationality of classifying spaces}, available at  \texttt{http://www.math.ucla.edu/$\sim$merkurev/publicat.htm}

\bibitem{Pukhlikov} A. V. Pukhlikov, {\it Rationality problem and birational rigidity}, Cohomological and geometric approaches to rationality problems, Progr. Math., {\bf 282}, Birkh\"auser, Boston, 2009, 275--311. 

\bibitem{roitman} A. Roitman, {\it Rational equivalence of zero-cycles}, Math. USSR Sbornik {\bf18} (1972) 571--588.

\bibitem{Rost} M. Rost, {\it Chow groups with coefficients}, Doc. Math. {\bf 1} (1996) no. 16, 319--393.


\bibitem{Saltman1} D. J. Saltman,
{\it Noether's problem over an algebraically closed  field}, Invent. math.
{\bf 77}
(1984), no. 1, 71-- 84.
\bibitem{Saltman} D. J. Saltman,
{\it Retract rational fields and cyclic Galois extensions,}  Israel J.
Math. {\bf 47} (1984), 165--215.

\bibitem{Srinivas} V. Srinivas, \emph{Gysin maps and cycle classes for Hodge cohomology}, Proc. Indian Acad. Sci. Math. Sci. {\bf 103} (1993), no.3, 209--247.

\bibitem{To15} B. Totaro, \emph{Hypersurfaces that are not stably rational},  arXiv:1502.04040, to appear in J. Amer. Math. Soc. 


\bibitem{Vial} Ch. Vial, {\it
Algebraic cycles and fibrations},  Documenta Math. {\bf  18} (2013) 1521--1553.

\bibitem{V13} C. Voisin,  \emph{Unirational threefolds with no universal codimension 2 cycles},
Invent.  Math. (1) {\bf 201} (2015) 207--237

\bibitem{VoisinRapport} C. Voisin, {\it Stable birational invariants and the L\"uroth problem}, available at \texttt{http://webusers.imj-prg.fr/$\sim$claire.voisin/Articlesweb/jdgvoisin.pdf}

\end{thebibliography}
\end{document}